\documentclass[12pt]{article}
\usepackage{hyperref}
\usepackage{cmap}                        
\usepackage[cp1251]{inputenc}            
\usepackage[english]{babel}
\usepackage[left=1in,right=1in,top=1in,bottom=1.5in]{geometry} 
\usepackage{amssymb,amsmath, amsthm, amscd,ifthen}
\usepackage{graphicx}
\usepackage{bm}

\newcommand{\N}{{\mathbb N}}
\newtheorem*{thm*}{Theorem}

\newcommand{\ff}{{\mathcal F}}
\newcommand{\G}{{\mathcal G}}

\newcommand{\aaa}{{\mathcal A}}

\newtheorem*{cla*}{Claim}
\newcommand{\bb}{{\mathcal B}}

\newtheorem{thm}{Theorem}

\newtheorem{lem}[thm]{Lemma}

\newtheorem{thrm}[thm]{Theorem}

\date{}

\title{Counting intersecting and pairs of cross-intersecting families}
\author{Peter Frankl, Andrey Kupavskii\footnote{Moscow Institute of Physics and Technology, Ecole Polytechnique F\'ed\'erale de Lausanne; Email: {\tt kupavskii@yandex.ru} \ \ Research of Andrey Kupavskii is supported by the grant RNF~16-11-10014.}}

\date{}

\begin{document}
\maketitle
\begin{abstract} A family of subsets of $\{1,\ldots,n\}$ is called {\it intersecting} if any two of its sets intersect. A classical result in extremal combinatorics due to Erd\H os, Ko, and Rado determines the maximum size of an intersecting family of $k$-subsets of $\{1,\ldots, n\}$. In this paper we study the following problem: how many intersecting families of $k$-subsets of $\{1,\ldots, n\}$ are there? Improving a result of Balogh, Das, Delcourt, Liu, and Sharifzadeh, we determine this quantity asymptotically for $n\ge 2k+2+2\sqrt{k\log k}$ and $k\to \infty$. Moreover, under the same assumptions we also determine asymptotically the number of {\it non-trivial} intersecting families, that is, intersecting families for which the intersection of all sets is empty. We obtain analogous results for pairs of cross-intersecting families.
\end{abstract}
MSc classification: 05D05

\section{Introduction}

A {\it family} is a collection of subsets of an $n$-element set $[n]$. Collections $\ff\subset{[n]\choose k}$ are called {\it $k$-uniform families}.
A family $\ff$ is called \textit{intersecting} if $F\cap F'\ne \emptyset$ holds for all $F,F'\in \ff$. Similarly, $\ff\subset {[n]\choose k}$ and $\G\subset{[n]\choose l}$ are called \textit{cross-intersecting}, if for all $F\in\ff$, $G\in \G$ one has $F\cap G\ne \emptyset.$

The research concerning intersecting families was initiated by Erd\H os, Ko and Rado, who determined the maximum size of intersecting families.

\begin{thrm}[Erd\H os, Ko, Rado \cite{EKR}] Suppose that $n\ge 2k>0$ and $\ff \subset{[n]\choose k}$ is intersecting. Then \begin{equation}\label{eqekr} |\ff|\le {n-1\choose k-1}.\end{equation}
\end{thrm}

The family of all $k$-sets containing a fixed element shows that (\ref{eqekr}) is best possible. Hilton and Milner proved in a stronger form that for $n>2k$ these are the only families on which the equality is attained. We say that an intersecting family $\ff$ is \textit{non-trivial} if $\bigcap_{F\in \ff} F=\emptyset,$ that is, if it cannot be pierced by a single point.

\begin{thrm}[Hilton, Milner \cite{HM}]\label{thmhm} Let $n>2k>0$ and suppose that $\ff\subset {[n]\choose k}$ is  a non-trivial intersecting family. Then \begin{equation}\label{eqhm} |\ff|\le {n-1\choose k-1}-{n-k-1\choose k-1}+1.\end{equation}
\end{thrm}

For $n = 2k+1$ the difference between the upper bounds \eqref{eqekr} and \eqref{eqhm} is only $k-1$. However, as $n-2k$ increases, this difference gets much larger. The number of subfamilies of $\ff$ is $2^{|\ff|}$, and thus the ratio between the number of subfamilies of the Erd\H os--Ko--Rado family and that of the Hilton--Milner family is $2^{k-1}$ for $n=2k+1$ and grows very fast as $n-2k$ increases. This serves as an indication that most intersecting families are \textit{trivial}, i.e., satisfy $\bigcap_{F\in \ff}F\ne \emptyset.$

In an important recent paper Balogh, Das, Delcourt, Liu, and Sharifzadeh \cite{BDD} proved this in the following quantitative form. Let $I(n,k)$ denote the \textit{total} number of intersecting families $\ff\subset {[n]\choose k}$.

\begin{thrm}[Balogh, Das, Delcourt, Liu, and Sharifzadeh \cite{BDD}]\label{thmbdd} If $n\ge 3k+8\log k$ then \begin{equation}\label{eqbdd} I(n,k) = (n+o(1))2^{{n-1\choose k-1}},\end{equation}
where $o(1)\to 0$ as $k\to \infty$.
\end{thrm}

One of the main tools of the proof of (\ref{eqbdd}) is a nice bound on the number of \textit{maximal} (i.e., non-extendable) intersecting families (see Lemma~\ref{lembb}). They obtain this bound using the following fundamental result of Bollob\'as.

\begin{thrm}[Bollob\'as \cite{Bo}] Suppose that $\aaa\subset {[n]\choose a}, \bb\subset{[n]\choose b}$ with $\aaa = \{A_1,\ldots, A_m\},\ \bb=\{B_1,\ldots, B_m\}$ satisfy $A_i\cap B_j=\emptyset$ iff $i=j$. Then \begin{equation}\label{eqbo} m\le {a+b\choose a}.\end{equation}
\end{thrm}

Note that the bound (\ref{eqbo}) is independent of $n$. In \cite{Bo} it is proved in a more general setting, not requiring uniformity, i.e., for $\aaa, \bb\subset 2^{[n]}$.
The uniform version (\ref{eqbo}) was rediscovered several years later by Jaeger--Payan \cite{JP} and Katona \cite{Ka2}.

We are going to use (\ref{eqbo}) to obtain an upper bound on the number of maximal pairs of cross-intersecting families.
Let us denote by $CI(n,a,b,t)$ ($CI(n,a,b,[t_1,t_2])$) the {\it number of pairs of cross-intersecting families} $\mathcal A\subset{[n]\choose a}, \mathcal B\subset {[n]\choose b}$ with $|\mathcal A|= t$ ($t_1\le |\aaa| \le t_2$).   We also denote $CI(n,a,b):=\sum_{t} CI(n,a,b,t)$.

We prove the following bound for the number of pairs of cross-intersecting families.
\begin{thm}\label{thmci} Choose $a,b,n\in \N$ and put $c:=\max\{a,b\}$, $T:={n-a+b-1\choose n-a}$. For  $n\ge a+b+2 \sqrt{c\log c}+2\max\{0,a-b\}$, $a,b\to \infty,$ and $b \gg \log a$ we  have
\begin{align}\label{eqci1}CI(n,a,b) =& (1+\delta_{ab}+o(1))2^{{n\choose c}},\\
\label{eqci2} CI(n,a,b,[1,T]) =& (1+o(1)){n\choose a} 2^{{n\choose b}-{n-a\choose b}},\end{align}
where $\delta_{ab}=1$ if $a=b$, and $0$ otherwise.
\end{thm}

In fact, \eqref{eqci2} easily implies \eqref{eqci1} (see Section~\ref{sec2} for details).\vskip+0.1cm

For a family $\ff\subset {[n]\choose k}$ we define the \textit{diversity} $\gamma(\ff)$ of $\ff$ to be $|\ff|-\Delta(\ff)$, where $\Delta(\ff):=\max_{i\in[n]}\big|\{F:i\in F\in \ff\}\big|$.
For an integer $t$ denote by $I(n,k, t)$ ($I(n,k,\ge t)$) the number of  intersecting families with diversity $t$ (at least $t$). In particular, $I(n,k,\ge 1)$ is the number of non-trivial intersecting families.
With the help of (\ref{eqci2}) we obtain a refinement of Theorem~\ref{thmbdd}.

\begin{thm}\label{thm1} For $n\ge 2k+2+2\sqrt{k\log k}$ and $k\to \infty$ we have
\begin{align}\label{eqi1} I(n,k) =& (n+o(1))2^{{n-1\choose k-1}},\\
\label{eqi2} I(n,k,\ge 1) =& (1+o(1))n{n-1\choose k} 2^{{n-1\choose k-1}-{n-k-1\choose k-1}}.\end{align}
\end{thm}

Again, it is easy to see that (\ref{eqi2}) implies (\ref{eqi1}) (see Section~\ref{sec3} for details).
In the next section we present the proof of Theorem~\ref{thmci}, and in Section~\ref{sec3} we give the proof of Theorem~\ref{thm1}.

\section{Cross-intersecting families}\label{sec2}
 Let us define the \textit{lexicographic} order on the $k$-subsets of $[n]$. We have $F\prec G$ in the lexicographic order if $\min F\setminus G<\min G\setminus F$ holds. E.g., $\{1,10\}\prec\{2,3\}$. For $0\le m\le {n\choose k}$ let $\mathcal L^{(k)}(m)$ denote the family of first $m$ $k$-sets in the lexicographic order. E.g., $\mathcal L^{(k)}\bigl({n-1\choose k-1}\bigr) = \{F\in{[n]\choose k}: 1\in F\}.$

Next we state the Kruskal--Katona  Theorem \cite{Kr}, \cite{Ka}, which is one of the most important results in extremal set theory.

\begin{thrm}[Kruskal \cite{Kr}, Katona \cite{Ka}]\label{thma} If $\mathcal A\subset{[n]\choose a}$ and $\bb\subset{[n]\choose b}$ are cross-intersecting then $\mathcal L^{(a)}(|\mathcal A|)$ and $\mathcal L^{(b)}(|\mathcal B|)$ are cross-intersecting as well.
\end{thrm}

Computationwise, the bounds arising from the Kruskal--Katona Theorem are not easy to handle. Lov\'asz \cite{L} found a slightly weaker but very handy form, which may be stated as follows.
\begin{thrm}[Lov\'asz \cite{L}]\label{thmb} Let $n\ge a+b$, and consider a pair of cross-intersecting families $\mathcal A\subset{[n]\choose a}$, $\bb\subset{[n]\choose b}$. If $|\mathcal A| =  {x\choose n-a}$ for a real number $x\ge n-a$, then
\begin{equation}\label{eq7} |\mathcal B|\le {n\choose b}- {x\choose b} \ \ \ \ \ \ \ \ \ \text{holds.}
\end{equation}
\end{thrm}

Note that for $x\ge k-1$ the polynomial ${x\choose k}$ is a monotone increasing function of $x$. Thus $x$ is uniquely determined by $|\mathcal A|$ and $a$.

We would also need the following result, which proof is based on Theorem~\ref{thma} and which is a combination of a result of Frankl and Tokushige \cite{FT} (Theorem 2 in \cite{FT}) and two results of Kupavskii and Zaharov \cite{KZ} (Part 1 of Theorem 1 and Corollary 1).
\begin{thrm}[Frankl, Tokushige, \cite{FT}, Kupavskii, Zakharov \cite{KZ}]\label{lemft} Let $n > a+b$ and suppose that the families $\mathcal A\subset{[n]\choose a},\mathcal B\subset{[n]\choose b}$ are cross-intersecting. If for some real number $\alpha \ge 1$ we have ${n-\alpha \choose n-a}\le |\mathcal A|\le {n-a+b-1\choose n-a}$, then \begin{equation}\label{eqft}|\mathcal A|+|\mathcal B|\le {n\choose b}+{n-\alpha \choose a-\alpha}-{n-\alpha \choose b}.\end{equation}
\end{thrm}

Note that the upper bound on $|\mathcal A|$ in this theorem is exactly the same as in \eqref{eqci2}.
\vskip+0.1cm
We go on to the proof of Theorem \ref{thmci}.
First we show that (\ref{eqci2}) implies (\ref{eqci1}). We may w.l.o.g. assume for this paragraph that $c=b\ge a$. For $b>a$ we have $T\ge {n\choose a}$ and thus $CI(n,a,b,t)=0$ for $t>T$. Therefore, we have $CI(n,a,b) = CI(n,a,b,0)+CI(n,a,b,[1,T])$. If $a=b$, then $T={n-1\choose a-1}$, and it follows from Theorem~\ref{thma} that if $\mathcal A,\mathcal B\subset{[n]\choose a}$ are cross-intersecting, then $\min\{|\aaa|,|\bb|\}\le {n-1\choose a-1}$. Therefore, in the case $a=b$ we have
$$2CI(n,a,b,0)-1 \le CI(n,a,b) \le 2(CI(n,a,b,0)+CI(n,a,b,[1,T])).$$ The ``-1'' in the first inequality stands for a pair of empty families, which is counted twice.   At the same time, we have $CI(n,a,b,0) = 2^{n\choose b}$. Thus, in both cases $b>a$ and $b=a$ it is sufficient to show that the right-hand side of (\ref{eqci2}) is $o(2^{n\choose b})$. We first note that $n-a-b\ge \sqrt n$ for $b\ge a$, $n\ge a+b+2\sqrt{b\log b}$, since $4b\log b\ge a+b+2\sqrt{b\log b}$. The rest is done by a simple calculation:

$$\frac{CI(n,a,b,[1,T])}{2^{{n\choose b}}}={n\choose a}2^{-{n-a\choose b}}\le 2^{n-{b+\sqrt{n}\choose b}} = o(1).$$
\vskip+0.1cm

Next, discuss the proof of the lower bound in (\ref{eqci2}). To obtain that many pairs of intersecting families, take $\mathcal A := \{A\}, A\in {[n]\choose a}$ and $\mathcal B(A):=\{B\in {[n]\choose b}: B\cap A\ne \emptyset\}$. Next, choose an arbitrary subfamily $\mathcal B\subset \bb(A)$. We only need to assure that few of these pairs of subfamilies are counted twice. Actually, we  count a pair of families  twice only in the case when $a=b$ and both $\mathcal A, \mathcal B$ consist of one set. The number of such pairs is ${[n]\choose a}^2$ and is negligible compared to the right hand side of \eqref{eqci2}.\\

We pass to the proof of the upper bound. \vskip+0.1cm

\pmb{$2\le |\mathcal A|\le n-a$}. Applying Theorem~\ref{thma}, the size of the (unique) maximal family $\mathcal B'$ that forms a cross-intersecting pair with $\aaa$ is maximized if $\mathcal A$ consists of two sets $A_1,A_2$ that intersect in $a-1$ elements. Therefore, $|\mathcal B'|\le {n\choose b}-{n-a+1\choose b}+{n-a-1\choose b-2}$. Any other family $\bb$ that forms  a cross-intersecting pair with $\mathcal A$ must be a subfamily of $\mathcal B'$.

So we can bound the number of pairs of cross-intersecting families $\mathcal A,\mathcal B$ with $2\le |\mathcal A|\le n-a$ as follows:
$$\frac{\sum_{t=2}^{n-a}CI(n,a,b,t)}{2^{{n\choose b}-{n-a\choose b}}}\le \sum_{t=2}^{n-a}{{n\choose a}\choose t}\frac{2^{{n\choose b}-{n-a+1\choose b}+{n-a-1\choose b-2}}}{2^{{n\choose b}-{n-a\choose b}}}\le 2^{n^2}2^{-{n-a-1\choose b-1}}=o(1).$$\vskip+0.1cm

\pmb{$n-a+1\le |\mathcal A|\le {n-u\choose n-a}$}, where $u=\sqrt{c\log c}+\max\{0,a-b\}$. Note that $n-a+1 = {n-a+1\choose n-a}$. In this case the bound is similar, but we use Theorem~\ref{thmb} to bound the size of $|\mathcal B|$. For $\mathcal A$ with $|\aaa|= t:={n-u'\choose n-a}$, where $u\le u'\le a-1$, we get $|\mathcal B|\le 2^{{n\choose b}-{n-u'\choose b}}$, and since ${{n\choose a}\choose t}\le 2^nt$, we have the following bound:

$$\frac{CI(n,a,b,t)}{2^{{n\choose b}-{n-a\choose b}}}\le {{n\choose a}\choose t}\frac{2^{{n\choose b}-{n-u'\choose b}}}{2^{{n\choose b}-{n-a\choose b}}}\le 2^{n{n-u'\choose n-a}}2^{-{n-u'-1\choose b-1}}.$$
At the same time we have $n\ge a+b+2u$ and
\begin{multline}\label{eq055}\frac{{n-u'\choose n-a}}{{n-u'-1\choose b-1}} = \frac {n-u'}{b}\prod_{i=0}^{n-a-b-1}\frac{n-b-u'-i}{n-a-i}\le \\ \frac nb \prod_{i=0}^{n-a-b-1}\frac{n-a-\sqrt{c\log c}-i}{n-a-i}\le \frac nb e^{-\sqrt{c\log c}(\sum_{i=b+1}^{n-a}\frac 1i)}\le \frac 1{2n}\end{multline}
for sufficiently large $c$. Indeed, $\sqrt{c\log c}\sum_{i=b+1}^{n-a}\frac 1i \ge \sqrt{c\log c}\sum_{i=b+1}^{b+2\sqrt{c\log c}}\frac 1i \ge (1+o(1))\frac {2(\sqrt{c\log c})^2}c=(2+o(1))\log c$, which justifies (\ref{eq055}) for $n\le b^{3/2}$. For $n>b^{3/2}$ we have $\sqrt{c\log c}\sum_{i=b+1}^{n-a}\frac 1i \ge (1+o(1))\sqrt{c\log c}\log \frac nb\ge (1+o(1))\sqrt{c\log c}\log n^{2/3}\gg \log n$, which justifies (\ref{eq055}) for $n>b^{3/2}$.\\

We conclude that $$\sum_{t=n-a+1}^{{n-u\choose n-a}}\frac{CI(n,a,b,t)}{2^{{n\choose b}-{n-b\choose b}}}\le 2^{-\frac 12{n-u'-1\choose b-1}} = o(1).$$
\vskip+0.1cm


\pmb{${n-u\choose n-a}< |\mathcal A|\le T$}, where $u=\sqrt{b\log b}+\max\{0,a-b\}$. Using the Bollobas set-pair inequality, it is not difficult to obtain the following bound on the number of maximal pairs of cross-intersecting families.

\begin{lem}\label{lembb} The number of maximal cross-intersecting pairs $\aaa'\subset{[n]\choose a},\bb'\subset{[n]\choose b}$ is at most  $[{n\choose a}{n\choose b}]^{{a+b\choose a}}$.
 \end{lem}
 We note that the proof is very similar to the proof of an analogous statement for intersecting families from \cite{BDD}.
\begin{proof} Find a minimal {\it $\mathcal B'$-generating} family $\mathcal M\subset\mathcal \aaa'$ such that $\mathcal \bb'=\{B\in {[n]\choose b}: B\cap M\ne \emptyset \text{ for all } M\in \mathcal M\}.$ We claim that $|\mathcal M|\le {a+b\choose a}$. Indeed, due to minimality, for each set $M'\in \mathcal M$ the family $\bb'':=\{B\in {[n]\choose b}: B\cap M\ne \emptyset \text{ for all }  M\in \mathcal M-\{M'\}\}$ strictly contains $\bb'$. Therefore, there is a set $B$ in $\bb''\setminus \bb'$ such that $B\cap M'=\emptyset, B\cap M\ne \emptyset $ for all $M\in \mathcal M-\{M'\}$. Applying the inequality (\ref{eqbo}) to $\mathcal M$ and the collection of such sets $B$, we get that $|\mathcal M|\le {a+b\choose b}$.

Interchanging the roles of $\aaa'$ and $\bb'$, we get that a minimal $\aaa'$-generating family has size at most ${a+b\choose b}$ as well. Now the bound stated in the lemma is just a crude upper bound on the number of ways one can choose these two generating families out of ${[n]\choose a}$ and ${[n]\choose b}$, respectively.
\end{proof}

Combined with the bound (\ref{eqft}) on the size of any maximal pair of families with such cardinalities, we get that
\begin{equation}\label{eq033}\frac{CI(n,a,b,[{n-u\choose n-a},T])}{2^{{n\choose b}-{n-a\choose b}}}\le
\bigg[{n\choose a}{n\choose b}\bigg]^{{a+b\choose a}}\frac{2^{{n\choose b}+{n-u\choose n-a}-{n-u\choose b}}}{2^{{n\choose b}-{n-a\choose b}}}\le 2^{2n{a+b\choose b}}2^{{n-u\choose n-a}-{n-u-1\choose b-1}}.\end{equation}


 We also have
$$\frac{{a+b\choose b}}{{n-u-1\choose b-1}} = \frac{n-u}{b}\prod_{i=0}^{b-1}\frac{a+b-i}{n-u-i}\le n\Big(\frac {a+b}{n-u}\Big)^{b}\le \frac 1{4n}.$$
Indeed, the last inequality is clearly valid for $n\ge (a+b)^2, b\to\infty$. If $n<(a+b)^2$, then the before-last expression is at most $$e^{\log n-\frac{b(n-a-b-u)}{n-u}}\le e^{2\log(a+b)-\frac{b(u+\max\{0,a-b\})}{O(a+b)}} \le e^{2\log(a+b)-\Omega(\min\{b,u\})}.$$
Since by the assumption we have $b\gg \log (a+b)$ and also, obviously, $u\gg \log (a+b)$, the last expression is at most $e^{-4\log(a+b)}<\frac 1{4n}.$\\

Taking into account (\ref{eq055}), which is valid for $u'=u$, we conclude that the right-hand side of (\ref{eq033}) is $o(1)$.

\section{Intersecting families}\label{sec3}
We need a theorem due to Frankl \cite{Fra1}, proved in the following, slightly stronger, form in \cite{KZ}.
\begin{thrm}[\cite{Fra1, KZ}]\label{thmc} Let $\ff\subset {[n]\choose k}$ be an intersecting family, and $n>2k$. Then, if $\gamma(\ff)\ge {n-u-1\choose k-u}$ for some real $3\le u\le k$, then \begin{equation}\label{eq01}|\ff|\le {n-1\choose k-1}+{n-u-1\choose k-u}-{n-u-1\choose k-1}.\end{equation}
\end{thrm}\vskip+0.2cm

We go on to the proof Theorem \ref{thm1}.  Let us first show that (\ref{eqi2}) implies (\ref{eqi1}). Indeed, using that $n-2k-1\ge \sqrt n$ and $k\to \infty$ in the assumptions of Theorem~\ref{thm1}, we get $$\frac{I(n,k,\ge 1)}{2^{n-1\choose k-1}}\sim n{n-1\choose k}2^{-{n-k-1\choose k-1}}\le 2^{n+\log n-{k+\sqrt n \choose k-1}}=o(1).$$
Therefore, $I(n,k,\ge 1) = o(2^{{n-1\choose k-1}})$. On the other hand, it is easy to see that $I(n,k,0) = (n+o(1))2^{{n-1\choose k-1}}$ (for the proof see \cite{BDD}).\\

Let us prove the lower bound in \eqref{eqi2}. For $S\in {[n]\choose k}, i\in [n]\setminus S$ define the family $\mathcal H(i,S): = \{S\}\cup \{H\in {[n]\choose k}: i\in H, H\cap S\ne \emptyset\}.$ Due to Theorem~\ref{thmhm}, these families are the largest non-trivial intersecting families. We have $|\mathcal H(i,s)| = {n-1\choose k-1}-{n-k-1\choose k-1}+1,$ and each such family contains no less than
\begin{equation}\label{eq003} 2^{{n-1\choose k-1}-{n-k-1\choose k-1}}-k 2^{{n-2\choose k-2}} = (1+o(1))2^{{n-1\choose k-1}-{n-k-1\choose k-1}}\end{equation}
non-trivial intersecting subfamilies, as $k\to \infty$. Indeed, a subfamily of $\mathcal H(i,S)$ containing $S$ is non-trivial unless all sets containing $i$ contain also a fixed $j\in S$. In other words, they must be a subset of a family $\mathcal I(i,j,S):=\{S\}\cup \{I\in{[n]\choose k}:i,j\in S\}$. The number of subfamilies of $\mathcal I(i,j,S)$ containing $S$ is $2^{{n-2\choose k-2}}$. Next, we have ${n-1\choose k-1}-{n-k-1\choose k-1}-{n-2\choose k-2}\ge {n-3\choose k-2}$, and thus the last inequality in the displayed formula above holds since $2^{{n-3\choose k-2}}\gg k$. Denote the set of all non-trivial subfamilies of $\mathcal H(i,S)$ by $\tilde{\mathcal H}(i,S)$.

Therefore, $\sum_{S\in {n\choose k}, i\notin S}|\tilde{\mathcal H}(i,S)| = (1+o(1))n{n-1\choose k}2^{{n-1\choose k-1}-{n-k-1\choose k-1}}.$ On the other hand, the pairwise intersections of these families are small: the families from $\tilde{\mathcal H}(i,S)\cap \tilde{\mathcal H}(i,S')$ form the set $I(n,k,2)$, and we do (somewhat implicitly) show in the proof that $I(n,k,2) = o(I(n,k,1))$. It could also be verified by a simple direct, but somewhat tedious calculation. Therefore, the lower bound is justified.\\



Next we prove the upper bound. For $i\in[n]$ and $\ff\subset {[n]\choose k}$ we use the standard notation \begin{align*}\ff(i):=&\{F-\{i\}:i\in F\in \ff\}\subset{[n]-\{i\}\choose k-1},\\
\ff(\bar i):=&\{F\in \ff:i\notin F\}\subset{[n]-\{i\}\choose k}.
\end{align*}
Note that if $\ff$ is intersecting then $\ff(i)$ and $\ff(\bar i)$ are cross-intersecting.

We count the number of families with different diversity separately. The number of families $\ff$ with $i$ being the most popular element and $\gamma(\ff)\le {n-3\choose k-2}$ is at most the number of cross-intersecting pairs $\ff(\bar i),\ \ff(i)$.

Therefore, we may apply (\ref{eqci2}) with $n':=n-1,a:=k, b:=k-1$, and get that the number of such families $\ff$ is at most $(1+o(1)){n-1\choose k} 2^{{n-1\choose k-1}-{n-k-1\choose k-1}}$. Note that $n'\ge a+b+2\sqrt{a\log a}+2$ and, in terms of Theorem \ref{thmci}, we have $T = {n-3\choose k-2}$ for our case. Multiplying the number of such families by the number of choices of $i$, we get the claimed asymptotic.

We are only left to prove that there are few families with diversity larger than ${n-3\choose k-2}$. Using the upper bound ${n\choose k}^{{2k-1\choose k-1}}$ for the number of maximal intersecting families in ${[n]\choose k}$ obtained in \cite{BDD} (see Lemma~\ref{lembb} for the proof of a similar statement), combined with the bound (\ref{eq01}) on the size of any maximal family with such diversity, we get that
\begin{equation}\label{eq03}\frac{I(n,k,\ge {n-3\choose k-2})}{2^{{n-1\choose k-1}-{n-k-1\choose k-1}}}\le {n\choose k}^{{2k-1\choose k-1}}\frac{2^{{n-1\choose k-1}+{n-4\choose k-3}-{n-4\choose k-1}}}{2^{{n-1\choose k-1}-{n-k-1\choose k-1}}}\le 2^{n{2k-1\choose k-1}}2^{{n-4\choose k-3}-{n-5\choose k-2}}.\end{equation}
Putting $n=2k+x$, we have ${n-4\choose k-3}/{n-5\choose k-2}= \frac{(n-4)(k-2)}{(n-k-1)(n-k-2)}\le\frac{(2k+x)k}{(k+x-2)^2} \le 1-\frac{x^2}{(k+x)^2}\le 1-   \frac 1k.$ On the other hand,
$$\frac{{2k-1\choose k-1}}{{n-5\choose k-2}} = \frac{n-4}{k-1}\prod_{i=1}^{k}\frac{2k-i}{n-3-i}\le n\Big(\frac {2k}{n-3}\Big)^{k}\le \frac 1{2kn},$$
where the last inequality is clearly valid for $n\ge 2k+2+2\sqrt{k\log k}$ and sufficiently large $k$. We conclude that the right-hand side of (\ref{eq03}) is at most $2^{\frac 1{2k}{n-5\choose k-2}}=o(1)$.

\end{document}